\documentclass[12pt]{amsart}
\usepackage{amssymb}
\usepackage[margin=1in]{geometry}
\usepackage{hyperref}

\newcommand{\modd}[1]{\,(\textup{mod}\, {#1})}
\def\Re{\textup{Re}\,}

\newtheorem{theorem}{Theorem}
\newtheorem{lemma}{Lemma}
\newtheorem{step}{Step}
\theoremstyle{remark}
\newtheorem*{remark}{Remark}

\begin{document}

\parskip5pt
\parindent20pt
\baselineskip15pt

\title[Determination of the Group generated by the Ratios $(an+b)/(An+B)$]{Harmonic Analysis on the Positive Rationals.  Determination of the Group generated by the Ratios $(an+b)/(An+B)$}

\author{P. D. T. A. Elliott}
\address{Department of Mathematics, University of Colorado Boulder, Boulder, Colorado 80309-0395 USA}
\email{pdtae@euclid.colorado.edu}

\author{Jonathan Kish}
\address{Department of Applied Mathematics, University of Colorado Boulder, Boulder, Colorado 80309-0526 USA}
\email{jonathan.kish@colorado.edu}





%
%

\maketitle

A systematic study of the representation of positive rationals by products and quotients of a given sequence of positive rationals was begun in the first author's Springer volume \cite{Elliott1985}.

Typically, representations of the form
\[
r^v = \prod_{j=1}^k \left( \frac{an_j+b}{An_j+B} \right)^{\varepsilon_j},
\]
where $a>0$, $A>0$, $b$, $B$ are integers for which $\Delta = aB - Ab \ne 0$, $\varepsilon_j = \pm 1$ and the integers $n_j$ are required to exceed a given value, are codified by the quotient group $G = \mathbb{Q}^*/\Gamma$, with $\mathbb{Q}^*$ the multiplicative group of positive rationals, $\Gamma$ its subgroup generated by the available $(an+b)/(An+B)$; formally a free abelian group on the primes with infinitely many relations.

It was shown in that volume that this particular group is finitely generated.  Homomorphisms of $G$ into the additive group of reals enabled a set of generators for the associated free group and membership of the attached torsion group to be determined;  the existence of an integer $v$ could be guaranteed but not its minimal value.

In the present paper it is shown how to determine $G$ through its dual.  In particular, the torsion group of $G$ will be a subgroup of the multiplicative reduced residue classes $\modd{2(aA\Delta)^3}$, hence its homomorphic image.  In the terminology of \cite{Elliott1985}, Chapter 22, suitably interpreted $G$ is \emph{arithmic}.  This respectively addresses a problem and settles in the affirmative a 1984 conjecture of the first author, c.f. \cite{Elliott1985}, Chapter 23, Unsolved Problems 11, 12.


Of particular importance is a recent result of Terence Tao, \cite{tao2015logarithmicaveraged}, on logarithmically weighted correlations of multiplicative functions.

It is convenient to begin with the following waystation of independent interest.

\begin{theorem} \label{elliott_kish_2015_groups_thm_01}
Let integers $a>0$, $A>0$, $b$, $B$, satisfy $\Delta = aB - Ab \ne 0$.  Set $\delta = 2(aA\Delta)^3$.

If a completely multiplicative complex-valued function $g$ satisfies
\[
g\left( \frac{an+b}{An+B} \right) = c \ne 0
\]
on all but finitely many positive integers, $n$, then there is a Dirichlet character $\modd{\delta}$ with which $g$ coincides on all primes that do not divide $\delta$.
\end{theorem}

\begin{remark}
Without loss of generality $(a,b) = 1 = (A,B)$.
\end{remark}

The argument is in four steps.


\begin{step} \label{elliott_kish_2015_groups_step_01}
There is an integer $m$ for which $g(p)^m=1$ whenever $(p,aA\Delta)=1$.  Moreover, $c$ is a root of unity.
\end{step}

Choose a rational $r$ for which $g(r)c=1$.  Then $g$ has the value $1$ on all but finitely many of the fractions $r(an+b)(An+B)^{-1}$, $n > 0$.  Let $\Gamma_0$ be the subgroup of $\mathbb{Q}^*$ that they generate.

We may regard $g$ as a character on the group $G_0 = \mathbb{Q}^*/\Gamma_0$.

Let $f$ be a homomorphism of $G_0$ into the additive reals, i.e., a completely additive function that satisfies
\[
f(an+b) - f(An+B) = -f(r), \quad n > n_0.
\]

According to \cite{Elliott1985}, Chapter 13, there is a real $H$ so that $f(m) = H\log m$ whenever $(m, aA\Delta)=1$.

\emph{There are two cases.  In the first case} there is a residue class for which $an+b, An+B$ are both prime to $aA\Delta$.  In this case
\[
H\log \left( \frac{an+b}{An+B} \right) = -f(r)
\]
for infinitely many positive integers, untenable unless $H=0$.

Thus $f(p) = 0$ if $(p, aA\Delta)=1$, and $f(r) = 0$.

For each such prime $p$ there is a positive integer $v$, possibly depending upon $p$, and a representation
\[
p^v = \prod_{j=1}^k \left( \frac{r(an_j+b)}{An_j+B} \right)^{\varepsilon_j}, \quad \varepsilon_j = \pm 1, n_j > n_0.
\]

As a consequence $g(p)^v = 1$.

Likewise $g(r)^w = 1$ for some positive integer, guaranteeing $c$ to be a root of unity.

If we choose a residue class $s \modd{aA\Delta}$ for which $((as + b)(As+B), aA\Delta)=1$, and consider the ratios
\[
(a (aA\Delta n + s) + b)(A(aA\Delta n + s) + B)^{-1},
\]
then the argument of \cite{Elliott1985}, Chapter 4, shows that the $p \modd{\Gamma_0}$ for primes $(p, aA\Delta)=1$ are finitely generated, and by primes also coprime to $aA\Delta$.

The existence of an integer $m$ so that the $g(p)$ are all $m^{th}$ roots of unity is now clear.

\emph{In the second case} we can find a residue class for which $an+b$ is divisible exactly by 2, but by no other prime divisor of $aA\Delta$, whilst $An+B$ is coprime to $aA\Delta$, c.f. \cite{Elliott1985}, Chapter 3; also a residue class for which the roles of the parameters $a$, $b$, $A$, $B$ are reversed.

Once again
\[
f(2) - H\log 2 + H\log \left( \frac{an+b}{An+b} \right) = -f(r)
\]
so that $H=0$, $f(2) = -f(r)$.

Moreover,
\[
-f(2) + H\log 2 + H \log \left( \frac{am+b}{Am+b} \right) = -f(r),
\]
so that $H=0$, $f(2) = f(r)$.

Hence $f(2) = 0 = f(r)$.

Since we may choose $r$ to have the form $2r_1$, where the rational $r_1$ is comprised  only of primes not dividing $aA\Delta$, $c$ itself is again a root of unity.

The inductive argument of \cite{Elliott1985}, Chapter 4, proceeds, since the class $2 \modd{\Gamma_0}$ has torsion, and for some $m$, $g(p)^m=1$ on all primes not dividing $aA\Delta$.


\begin{step}\label{elliott_kish_2015_groups_step_02}
There is a Dirichlet character $\chi \modd{D}$ and a set of primes $q$ with $\sum q^{-1}$ convergent, such that $g(p) = \chi(p)$ on all remaining primes.
\end{step}

\begin{remark}
Without loss of generality we may assume $\chi$ to be primitive.
\end{remark}

The following is Theorem 1.3 of Tao, \cite{tao2015logarithmicaveraged}.


\begin{lemma} \label{elliott_kish_2015_groups_lem_A}
Let the integers $a>0$, $b>0$, $c$, $d$, satisfy $ad-bc \ne 0$.  Let $\varepsilon >0$ and suppose that $A_0$ is sufficiently large depending upon $\varepsilon$, $a$, $b$, $c$, $d$.  Let $x \ge w \ge A_0$ and let $g_1$, $g_2$ be multiplicative functions, with values in the complex unit disc, for which
\[
\sum_{p \le x} p^{-1} (1 - \Re g_1(p) \overline{\chi}(p) p^{-it}) \ge A_0
\]
for all Dirichlet characters of period at most $A_0$, and all real numbers $t$ with $|t| \le A_0x$.

Then
\[
\left| \sum_{x/w < n \le x} n^{-1} g_1(a n + b) g_2(c n + d) \right| \le \varepsilon \log w.
\]
\end{lemma}

To implement the result of Tao we apply


\begin{lemma} \label{elliott_kish_2015_groups_lem_B}
If $|g(p)| \le 1$ on the primes and
\[
\sum_{p \le x} p^{-1} ( 1 - \Re g(p) p^{i\lambda(x)}) \ll 1,
\]
with $\lambda(x)$ real and $\lambda(x) \ll x$, for $x \ge 2$, then there is a constant $\alpha$ for which $\lambda(x) = \alpha + O((\log x)^{-1})$ and the series $\sum p^{-1} (1- \Re g(p) p^{i\alpha})$ converges.
\end{lemma}


A somewhat elaborate version of Lemma \ref{elliott_kish_2015_groups_lem_B} is employed in the first author's study of correlations attached to the sums of renormalised shifted additive functions, \cite{elliott1994correlationmemoir}, Section 10.3.  The present version may be found as Lemma 17 of \cite{elliott2010valuedistribution}.  Note that on page 84 line 2 of that account the first sum over the primes contains a surplus factor of $g(p)$.

Applying Lemma \ref{elliott_kish_2015_groups_lem_A} to $g(an+b) \overline{g}(An+B)$ guarantees a constant $A_0$ and for each $x$ sufficiently large a pair $\chi \modd{D_x}$, $t_x$ real, with $\chi$ a Dirichlet character to a modulus $D_x \le A_0$, $|t_x| \le A_0x$, for which the sums
\[
\sum_{p \le x} p^{-1} (1 - \Re g(p) \overline{\chi}(p) p^{-it})
\]
are uniformly bounded.

The characters belong to a finite set.  This provides a positive integer $k$ for which $\overline{\chi}(p)^k=1$ on all but finitely many primes.  The inequality $1-\Re z^k \le k^2(1-\Re z)$, valid in the complex unit disc, shows the sums
\[
\sum_{p \le x} p^{-1} ( 1 - \Re g(p)^k p^{-ikt})
\]
to be uniformly bounded.

An application of Lemma \ref{elliott_kish_2015_groups_lem_B} guarantees a real $\beta$ for which $t_x = \beta + O((\log x)^{-1})$.  Hence
\[
\sum_{p \le x} ( 1 - \Re g(p) \overline{\chi}(p)p^{-it_x})
\]
differs from a similar sum with $t_x$ replaced by $\beta$, by 
\[
\ll \sum_{p \le x} p^{-1} |p^{-it_x} - p^{-i\beta}| \ll \sum_{p \le x} p^{-1} |t_x - \beta| \log p \ll 1.
\]

The sums
\[
\sum_{p \le x} p^{-1} ( 1 - \Re g(p) \overline{\chi}(p)p^{-i\beta})
\]
are uniformly bounded and, since there are only finitely many possibilities for the character, for some character the corresponding infinite series converges.

This result is interesting in its own right, requiring only that
\[
\liminf_{x \to \infty} (\log x)^{-1} \left| \sum_{n \le x} n^{-1} g_1(a n + b) g_2(A n + B) \right| > 0.
\]

To complete Step \ref{elliott_kish_2015_groups_step_02} we note that by Step \ref{elliott_kish_2015_groups_step_01}, $g(p)^m=1$ on all but finitely many primes.  Hence the series
\[
\sum p^{-1} ( 1 - \Re p^{-imk\beta})
\]
converges; and that is only tenable if $\beta=0$.

On the primes for which $g(p) \ne \chi(p)$, $g(p) \overline{\chi}(p)$ is a nontrivial $mk^{th}$ root of unity and the corresponding summands are uniformly bounded from below.

%


\begin{step}\label{elliott_kish_2015_groups_step_03}
$g(p) = \chi(p)$ provided $(p,D)=1$, $(p,aA\Delta)=1$.
\end{step}

It is convenient to establish a more general result.

\begin{lemma} \label{elliott_kish_2015_groups_lem_01}
Let the integers $a>0$, $A>0$, $b$, $B$ satisfy $\Delta = aB - Ab \ne 0$.  Let $g$ be a completely multiplicative function, with values in the complex unit disc, that satisfies 
\[
\limsup_{x \to \infty} x^{-1} \left| \sum_{n \le x} g(an+b) \overline{g}(An+B) \right| = 1,
\]
it being understood that finitely many of the summands may be omitted.

If $g$ is 1 on all primes that do not belong to a set of primes $q$ for which $\sum q^{-1}$ converges, then $g$ is 1 on all the primes that do not divide $\bigl( \frac{a}{(a,b)}, \frac{A}{(A,B)}\bigr)$.
\end{lemma}

\begin{remark}
Without loss of generality we may assume $(a,b)=1=(A,B)$.
\end{remark}

The next result is a particular case of a theorem of Stepanauskas, \cite{stepanasukas2002meanvaluesV}, who was concerned with allowing the parameters $a_j$, $b_j$ and functions $g_j$, $j = 1, 2$, to grow with the variable, $x$.


\begin{lemma} \label{elliott_kish_2015_groups_lem_C}
Let $g_1$, $g_2$ be multiplicative arithmetic functions with values in the complex unit disc.  Define the multiplicative functions $h_j$ by $h_j(p^m) = g_j(p^m) - g_j(p^{m-1})$, $m = 1, 2, \dots$, $j = 1, 2$.  Let $a_1$, $a_2$, $b_1$, $b_2$ be integers satisfying $a_1 > 0$, $a_2 > 0$, $(a_j,b_j)=1$, $j = 1, 2$, $\Delta = a_1b_2 - a_2b_1 \ne 0$.

Define
\[
w_p = \mathop{ \sum_{m_1=0}^\infty \sum_{m_2=0}^\infty}_{\substack{(p_j, a_j)=1, j = 1, 2 \\ (p^{m_1}, p^{m_2}) \mid \Delta}} \frac{h_1(p^{m_1}) h_2(p^{m_2})}{[p^{m_1}, p^{m_2}]},
\]
\[
S(x) = \sum_{j=1}^2 \sum_{\log x < p \le x} \frac{|s_j(p)|^2}{p},
\]
where
\[
s_j(p) = \begin{cases} g_j(p)-1 & \text{if} \ p \nmid \Delta a_j, \\ g_1(p)g_2(p) -1 & \text{if} \ p \mid \Delta, p \nmid a_j, \\ 0 & \text{if} \ p \mid a_j. \end{cases}
\]

Then,
uniformly in $x \ge 2$,
\[
x^{-1} \sum_{n \le x} g_1(a_1 n + b_1) g_2(a_2 n+b_2) - \prod_{p \le x} w_p \ll S(x)^{1/2} + (\log x)^{-1},
\]
the implied constant independent of the $g_j$.
\end{lemma}

\noindent \emph{Proof of Lemma \ref{elliott_kish_2015_groups_lem_01}}.  We define the multiplicative function $h$ by Dirichlet convolution, $g = 1 * h$, so that
\[
h(p^m) = g(p^m) - g(p^{m-1}) = (g(p)-1) g(p)^{m-1}, \quad m = 1, 2, \dots
\]
Then by Lemma \ref{elliott_kish_2015_groups_lem_C},
\[
\lim_{x \to \infty} x^{-1} \sum_{n \le x} g(an+b) \overline{g}(An+B) = \prod_q w_q,
\]
where 
\[
w_q = 1 + q^{-\Delta_q} 2\Re ((g(q)-1)(q-g(q))^{-1})
\]
if $q \nmid aA$, $\Delta_q$ being the highest power to which $q$ divides $\Delta$.  Otherwise $w_q = (q-1)(q-g(q))^{-1}$ if $q \nmid a$, $q \mid A$; and its conjugate if $q \mid a$, $q \nmid A$.

Since this result holds for any choices of the $g(q)$, we may set all but one $g(q)=0$ and conclude that $|w_q| \le 1$.  The displayed hypothesis of Lemma \ref{elliott_kish_2015_groups_lem_01} then shows that $|w_q| = 1$ for each prime $q$.

If $(q, aA)=1$, then $w_q$ is real, hence $\pm 1$.  This requires
\[
\Re \left( \frac{1-\frac{1}{q}}{1-\frac{g(q)}{q}} \right) = 1- q^{\Delta_q} \quad \text{or} \quad 1.
\]
As a diagram in the complex plane shows, the left hand ratio is in absolute value less than 1 unless $g(q)=1$, and is certainly not zero.

The remaining cases are similar but simpler.

\noindent \emph{Completion of Step \ref{elliott_kish_2015_groups_step_03}}.  The argument differs slightly according to the circumstances of the cases considered in Step \ref{elliott_kish_2015_groups_step_01}.

In the first case we choose a residue class $s \modd{aA\Delta}$ for which $((as+b)(As+B), aA\Delta)=1$ and apply Lemma \ref{elliott_kish_2015_groups_lem_01} to the function
\[
g \overline{\chi} (a(aA\Delta n+s) + b) \overline{g}\chi(A(aA\Delta n + s) + B).
\]
On the primes for which $(p, DaA\Delta)=1$, $g\overline{\chi}=1$.

In the second case(s) we adopt the modifications employed in Step \ref{elliott_kish_2015_groups_step_01}, noting that $|g(2)|=1$.


\begin{step}\label{elliott_kish_2015_groups_step_04}
$D$ divides $2(aA\Delta)^3$.
\end{step}

To this end we apply the following result.

\begin{lemma} \label{elliott_kish_2015_groups_lem_02}
%
Let the integers $u_j >0$, $v_j$, $(u_j, v_j)=1$, $j=1,2$ satisfy $\Delta_1 = u_1v_2 - u_2v_1 \ne 0$.  Assume that the primitive Dirichlet character $\chi_D$ satisfies
\[
\chi_D \left( \frac{u_1k + v_1}{u_2k+v_2} \right) = c \ne 0
\]
for all $k$ such that $(u_jk+v_j, D) = 1$, $j=1, 2$, and that there exists a $k_0$ for which this holds;  hence a class $k_0 \modd{D}$.  

Then $D \mid (u_1, u_2)\Delta_1$ if $D$ is odd, $D \mid 2(u_1, u_2)\Delta_1$ otherwise.
\end{lemma}

\noindent \emph{Proof of Lemma \ref{elliott_kish_2015_groups_lem_02}}.  Define
\[
\chi_D \left( \frac{u_1k + v_1}{u_2k+v_2} \right) = 0 \quad \text{if} \quad ((u_1k+v_1)(u_2k+v_2),D) > 1.
\]

If $D = \prod_{p^t \| D} p^t$, then there is a decomposition $\chi_D = \prod_{p^t \| D} \chi_{p^t}$.

Correspondingly
\[
\chi_D \left( \frac{u_1k + v_1}{u_2k+v_2} \right) = \chi_{p^t} \left( \frac{u_1k + v_1}{u_2k+v_2} \right) \chi_{D_1} \left( \frac{u_1k + v_1}{u_2k+v_2} \right)
\]
where $D_1 = p^{-t}D$.

If we set $k = \tilde{k} D_1 + k_0$, then $(u_j(\tilde{k}D_1 + k_0) + v_j, D_1) = 1$, $j=1, 2$, hence
\[
\chi_{D_1} \left( \frac{u_1(\tilde{k}D_1 + k_0) + v_1}{u_2(\tilde{k}D_1 + k_0) + v_2} \right) = \chi_{D_1} \left( \frac{u_1k_0 + v_1}{u_2k_0 + v_2} \right) \ne 0.
\]
Replacing $\tilde{k}$ by $k$:
\[
\chi_{p^t} \left( \frac{u_1(kD_1 + k_0) + v_1}{u_2(kD_1 + k_0) + v_2} \right) = \begin{cases} c_1 \ne 0 & \text{if} \ (u_j(kD_1 + k_0) + v_j,p) = 1, j = 1, 2 \\ 0 & \text{otherwise.} \end{cases}
\]

We have reduced ourselves to the case $D = p^t$, with $u_j, v_j$ replaced by $u_jD_1$, $k_0u_j+v_j$, $j = 1, 2$.

Note that $D_1u_1(k_0u_2 + v_2) - D_1 u_2(k_0 u_1 + v_1) = D_1\Delta_1$.

For convenience of exposition, write $w_j$ for $k_0 u_j + v_j$, $j=1, 2$.

Assume that $p^{s_j} \| u_j$ with, without loss of generality, $s_2 \le s_1$.  Otherwise, consider $\overline{\chi}_{p^t}$.

\emph{If $s_2 \ge t$ we have $p^t \mid (u_1, u_2)$ and we (temporarily) stop}.  Otherwise, set $u_j = p^{s_j} m_j$, so that $p \nmid m_j$, $j = 1, 2$.  Then
\begin{align*}
 \frac{u_1D_1k + w_1}{u_2D_1k+w_2} & = \frac{m_2(u_1D_1k + w_1)}{m_2(u_2D_1k+w_2)} \\
&  = \frac{m_1 p^{s_1 - s_2} p^{s_2} m_2 D_1 k + m_2 w_1}{m_2(m_2 p^{s_2}D_1k + w_2)} \\
& = \frac{m_1p^{s_1 - s_2}( m_2 p^{s_2} D_1 k+ w_2) + m_2 w_1 - m_1 p^{s_1 - s_2} w_2}{m_2(m_2 p^{s_2} D_1 k + w_2)} \\
& = \frac{m_1}{m_2} p^{s_1 - s_2} + \frac{m_2w_1 - m_1 p^{s_1 - s_2} w_2}{m_2(u_2 D_1 k + w_2)}.
\end{align*}
Since $\chi_{p^t}$ is primitive, for an appropriate Gauss sum $\varepsilon(\chi_{p^t})$ with $|\varepsilon(\chi_{p^t})| = p^{t/2}$,
\[
\chi_{p^t} \left( \frac{u_1D_1k + w_1}{u_2D_1k + w_2} \right) \varepsilon(\chi_{p^t}) = \sum\limits_{r=1}^{p^t} \overline{\chi}_{p^t}(r) \exp \left( \frac{2\pi i r}{p^t} \left( \frac{u_1 D_1k + w_1}{u_2D_1 k + w_2} \right) \right)
\]
whenever $(u_2D_1k + w_2, p) = 1$.

In particular,
\begin{equation}
\varepsilon(\chi_{p^t}) \sideset{}{'}\sum\limits_{k=1}^{p^{t-s_2}} \chi_{p^t} \left( \frac{u_1D_1k + w_1}{u_2D_1k + w_2} \right) = \sum\limits_{r=1}^{p^t} \overline{\chi}_{p^t}(r)  \sideset{}{'}\sum\limits_{k=1}^{p^{t-s_2}}  \exp \left( \frac{2\pi i r}{p^t} \left( \frac{u_1 D_1k + w_1}{u_2D_1 k + w_2} \right) \right), \label{elliott_kish_2015_groups_01_fund_relation}
\end{equation}
where $'$ denotes that summation is confined to terms with $(u_2D_1k + w_2,p)=1$.

The second innersum has the alternative representation
\[
M = \sideset{}{'}\sum\limits_{k=1}^{p^{t-s_2}} \exp \left( \frac{2\pi i r}{p^t} \left( \frac{m_1}{m_2} p^{s_1 - s_2} + \frac{L}{u_2 D_1 k + w_2} \right) \right)
\]
where
\[
L = \overline{m}_2(m_2w_1 -  m_1 p^{s_1 - s_2} w_2), \qquad m_2\overline{m}_2 \equiv 1 \modd{p^t}.
\]
Here $1/(u_2D_1k + w_2)$ is likewise interpreted as a group inverse $\modd{p^t}$.

\emph{For ease of notation we replace $s_2$ by $s$}.

Assume that $s \ge 1$.  The restriction $(u_2D_1k + w_2,p)=1$ is then automatically satisfied.

For $1 \le k \le p^{t-s}$, arguing via representations, we map the class $u_2D_1k + w_2 \modd{p^t}$ onto the class $j_k \modd{p^{t-s}}$ given by
\[
j_k = \frac{\overline{m_2p^sD_1k + w_2} - \overline{w}_2}{p^s},
\]
the inverses taken $\modd{p^t}$.  Since
\[
w_2( m_2 p^s D_1 k + w_2)(\overline{m_2p^s D_1 k + w_2} - \overline{w}_2) \equiv m_2p^s D_1 k \modd{p^t},
\]
$j_k$ is well-defined.  It is the class
\[
(w_2(m_2p^s D_1 k + w_2))^{-1} m_2 D_1 k \modd{p^{t-s}},
\]
the group inverse $^{-1}$ here taken in the reduced residue class group $\modd{p^{t-s}}$.

Moreover, if $j_{k_1} \equiv j_{k_2} \modd{p^{t-s}}$, $1 \le k_1 \le k_2 \le p^{t-s}$, then
\[
\overline{m_2p^s D_1 k_1 + w_2} - \overline{w}_2 \equiv \overline{m_2p^s D_1 k_2 + w_2} - \overline{w_2} \modd{p^t}, 
\]
from which $k_1 \equiv k_2 \modd{p^{t-s}}$ rapidly follows.  The map is one-to-one and covers every class $\modd{p^{t-s}}$.

In this case
\begin{align*}
M & = \sum\limits_{j=1}^{p^{t-s}} \exp\left( \frac{2\pi i r}{p^t} \left( \frac{m_1}{m_2} p^{s_1 - s} + L(p^sj + \overline{w}_2) \right) \right) \\
& = \exp\left( \frac{2\pi i r}{p^t} \left( \frac{m_1}{m_2}p^{s_1-s} + L\overline{w}_2 \right) \right) \sum\limits_{j=1}^{p^{t-s}} \exp\left( \frac{2\pi i r L j}{p^{t-s}} \right) = 0,
\end{align*}
unless $p^{t-s} \mid L$.  Note that from \eqref{elliott_kish_2015_groups_01_fund_relation} we may assume that $(r,p)=1$.

Hence $p^t \mid p^s L$, $p^t \mid (p^sm_2w_1 - p^{s_1}m_1w_2)$, i.e., $p^t \mid (u_1v_2 - v_1u_2)$, and once more we (temporarily) stop.

\emph{Variant}.  Suppose now that $s = s_2 = 0$, i.e., $p \nmid u_2$.  In this case
\[
M =  \sideset{}{'}\sum\limits_{k=1}^{p^t}  \exp\left( \frac{2\pi ir}{p^t} \left( \frac{u_1}{u_2} + \frac{L}{u_2D_1k + w_2} \right) \right).
\]
Set $\gamma = \exp(2\pi i r u_1 / p^t u_2)$.  Working within the reduced residue class group $\modd{p^t}$, we introduce a new variable $z = \overline{u_2D_1k + w_2}$.  Since $z \to \overline{z}$ permutes the group, $M$ has a representation
\[
M = \gamma \sum\limits_{\substack{z = 1 \\ (z,p)=1}}^{p^t} \exp\left( \frac{2\pi i rL z}{p^t} \right),
\]
a Ramanujan sum with an alternative representation in terms of the M\"obius function:
\[
M = \gamma \sum\limits_{d \mid (p^t, rL)} \mu\left( \frac{p^t}{d} \right) d.
\]

If $p^{t-1} \nmid L$, then $d = p^h$ with $h \le t-2$, $\mu(p^t d^{-1}) = 0$ and, from \eqref{elliott_kish_2015_groups_01_fund_relation},
\[
 \sideset{}{'}\sum\limits_{k=1}^{p^t} \chi_{p^t} \left( \frac{u_1D_1k + w_1}{u_2D_1k+w_2} \right) = 0.
\]

Since the summand with $k=p^t$ is nonzero, this is impossible.

%
%

Therefore $p^{t-1} \mid L$.

\emph{If $p^t \mid L$, then we stop, for once again $p^t \mid \overline{m}_2 (m_2w_1 - m_1w_2)$, $p^t \mid (u_2w_1 - u_1w_2)$, i.e., $p^t \mid (u_1v_2 - u_2v_1)$}.

Otherwise, $p^{t-1} \| L$, and the Ramanujan sum has value $-p^{t-1}$.

The fundamental relation \eqref{elliott_kish_2015_groups_01_fund_relation} becomes
\begin{align*}
\varepsilon(\chi_{p^t}) \sum\limits_{k=1}^{p^t} \chi_{p^t} \left( \frac{u_1D_1k + w_1}{u_2D_1k + w_2} \right) & = -p^{t-1} \sum\limits_{r=1}^{p^t} \overline{\chi}_{p^t}(r) \exp\left( \frac{2\pi i ru_1}{p^t u_2} \right) \\
& = -p^{t-1} \chi_{p^t} \left( \frac{u_1}{u_2} \right) \varepsilon(\chi_{p^t} );
\end{align*}
we may cancel the gaussian factors.

Several arguments now present themselves.  For example, in absolute value the left-hand sum is at least $p^t - (2p^{t-1}-1)$, guaranteeing that the prime $p$ is 2.

Otherwise, we note that $p \nmid u_1$ and set $y = \overline{D_1k+ \overline{u}_2 w_2}$ to obtain a representation
\[
\sum\limits_{\substack{y=1 \\ (y,p)=1}}^{p^t} \chi_{p^t} ( 1 + y(\overline{u}_1w_1 - \overline{u}_2w_2)) = -p^{t-1}\chi_{p^t}(u_1 \overline{u}_2).
\]

Hence $c(p^t - p^{t-1}) = -p^{t-1} \chi_{p^t}(u_1 \overline{u}_2)$, $p=2$ and $c = -\chi_{p^t}(u_1 \overline{u}_2)$.


At this stage, our initial hypothesis implies one of 
\begin{enumerate}
\item[(i)]  $p^t \mid (u_1, u_2)$,
\item[(ii)]  $p^t \mid (u_1v_2 - u_2v_1)$ if $p \ge 3$,
\item[(iii)]  $p^{t-1} \mid (u_1v_2 - v_1 u_2)$ if $p=2$.
%
%
\end{enumerate}

The conclusion of Lemma \ref{elliott_kish_2015_groups_lem_02} is now clear.


\noindent \emph{Completion of Step \ref{elliott_kish_2015_groups_step_04}}.  Once again there are small modifications according to the cases of Step \ref{elliott_kish_2015_groups_step_01}.  In the notation for the first case of Step \ref{elliott_kish_2015_groups_step_03},
\begin{align*}
1 & = g \left( \frac{a(aA\Delta n + s)+b}{A(aA\Delta n + s)+B} \right) \\
& = \chi \left( \frac{a(aA\Delta n + s)+b}{A(aA\Delta n + s)+B} \right), \quad \chi\modd{D},
\end{align*}
as long as $\chi$ is nonzero.  By Lemma \ref{elliott_kish_2015_groups_lem_02}, $D$ divides $2(a,A)(aA)^2 \Delta^3$.

The remaining cases proceed similarly.

Theorem \ref{elliott_kish_2015_groups_thm_01} is established.


\noindent \textbf{Constraints}.  Given integers $a>0$, $A>0$, $b$, $B$ with $\Delta = aB - Ab \ne 0$, define $\alpha = (a,b)$, $\beta = (A,B)$, $a_1 = a\alpha^{-1}$, $b_1 = b\alpha^{-1}$, $A_1 = A\beta^{-1}$, $B_1 = B\beta^{-1}$, $\rho = \alpha \beta^{-1}$ in lowers terms and $\Delta_1 = a_1 B_1 - A_1 b_1$, so that $\Delta = \alpha \beta \Delta_1$.  Define $\rho_0$ to be that part of $\rho$ made up of primes that divide $(a_1, A_1)$, to be 1 if there are none.

Assume that Theorem \ref{elliott_kish_2015_groups_thm_01} is in operation with $c=1$.

Since there are only two constraints upon $n$, we can find a positive $n$, and so a complete residue class, for which $(a_1 n + b_1)(A_1 n + B_1)$ is not divisible by a given odd prime, $p$.

Moreover, either this is possible for $p=2$ or there is a pair of integers $n_1$, $n_2$ for which $2 \| (a_1 n_1 + b_1)$, $\tfrac{1}{2} (a_1n_1 + b_1)(A_1 n_1 + B_1)$ is odd, and $2 \| (A_1n_2+B_1)$, $(a_1 n_2 + b_1)\tfrac{1}{2} (A_1 n_2 + B_1)$ is odd.

Hence we can either arrange for $(an+b)/(An+b) = \rho r$ with $(r,\delta)=1$, or for the pair $(an_1+b)/(An_1+B)= 2\rho r_1$, $(an_2+b)/(An_2+B) = (\rho/2)r_2$, with $(r_1r_2, \delta)=1$.

In either case, $g(p)^{2\phi(\delta)}=1$.

For each prime divisor $p$ of $\delta$ that does not divide $(a_1, A_1)$ we can arrange a value of $n$ such that $p \|  (a_1 n + b_1)/(A_1 n + B_1)$.  If $p^z \| \Delta_1$ and $p \nmid a_1$ (say), then we choose $n$ such that $p^{z+1} \| (a_1 n + b_1)$.  The identity $a_1(A_1 n + B_1) - A_1(a_1n + b_1) = \Delta_1$ shows that $p^z \| (A_1 n + B_1)$.  

Application of the Chinese Remainder Theorem together with the above argument shows that again $g(p)^{2\phi(\delta)}=1$.


\noindent \textbf{A characterisation of $G$}.  Let $g$ be a completely multiplicative function with values in the complex unit disc and that coincides with a Dirichlet character $\modd{\delta}$ on the integers prime to $\delta$.


Then
\[
\sum_{n \le x} g(an+b) \overline{g}(An+B) = g(\rho) \sum_{\substack{d_j \mid \delta_\infty \\ (d_1,d_2) \mid \Delta}} g(d_1) \overline{g}(d_2) \sideset{}{'}\sum_{n \le x} \chi \left( \frac{a_1 n+b_1}{d_1} \right) \overline{\chi} \left( \frac{A_1n+B_1}{d_2} \right)
\]
where $d_j \mid \delta_\infty$ denotes that $d_j$ is comprised of powers of the primes that divide $\delta$ and the innersum is taken over integers $n$ for which $a_1 n + b_1$ is divisible by $d_1$ and $A_1 n + B_1$ by $d_2$.  The value of the innermost summand is determined by the residue class $\modd{\delta[d_1, d_2]}$ to which $n$ belongs.

A typical innersum has the uniform bound $O( x/\max(d_1, d_2))$, and the asymptotic estimate
\[
\sideset{}{'}\sum_{n \modd{\delta[d_1,d_2]}} \chi\left( \frac{a_1 n + b_1}{d_1} \right) \overline{\chi}\left( \frac{A_1 n + B_1}{d_2} \right) \left( \frac{x}{\delta[d_1,d_2]} + O(1) \right).
\]
In particular, the sum
\[
\theta_{d_1, d_2}(\chi) = \frac{1}{\delta[d_1, d_2]} \sideset{}{'}\sum_{n \modd{\delta[d_1,d_2]}}  \chi\left( \frac{a_1 n + b_1}{d_1} \right) \overline{\chi}\left( \frac{A_1 n + B_1}{d_2} \right) 
\]
is $O((\max(d_1, d_2)^{-1})$.

Noting that
\[
\sum_{\substack{d \le y \\ d \mid \delta_\infty}} 1 \le \prod_{p \mid \delta} \left[ \frac{\log y}{\log p} \right] \ll (\log y)^\omega
\]
where $\omega$ denotes the number of distinct prime divisors of $\delta$, we see that
\[
\sum_{\sqrt{y} < [d_1, d_2] \le y} | \theta_{d_1, d_2}(\chi) | \ll (\log y)^\omega \sum_{\substack{d \mid \delta_\infty \\ d > y^{1/4}}} \frac{1}{d}
\]
\[
\ll (\log y)^\omega \sum_{p \mid \delta} \sum_{p^m > y^{1/4\omega}} \frac{1}{p^m} \ll (\log y)^\omega y^{-1/4\omega}.
\]
Hence
\[
\lim_{x \to \infty} x^{-1} \sum_{n \le x} g(an+b) \overline{g}(An+B) = g(\rho) \sum_{\substack{d_j \mid \delta_\infty \\ (d_1, d_2)\mid \Delta}} g(d_1) \overline{g}(d_2) \theta_{d_1, d_2}(\chi),
\]
which sum we shall denote by $S(g,\chi)$.

We have so far shown that for any character $g: G \to \{ z \in \mathbb{C}, |z|=1\}$,
\begin{enumerate}
\item[(i)]  \emph{There is a Dirichlet character $\chi\modd{\delta}$ so that $g(p) = \chi(p)$ whenever $(p,\delta)=1$},
\item[(ii)] \emph{ $g(p)^{2\phi(\delta)} = 1$ if $p \nmid (a_1, A_1)$; $g(\rho)^{2\phi(\delta)}=1$},
\item[(iii)]   $\displaystyle \mathop{\, g(\rho) \sum g(d_1)}_{\substack{d_j \mid \delta_\infty \\ (d_1, d_2) \mid \Delta}}\overline{g}(d_2) \theta_{d_1, d_2}(\chi) = 1 $
\end{enumerate}

Conversely, suppose that we choose a $\chi \modd{\delta}$ and set $g(p) = \chi(p)$ if $(p,\delta)=1$.  If we can, we next choose the $g(p)$ for $p$ dividing $\delta$ but not $(a_1, A_1)$ to satisfy $g(p)^{2\phi(\delta)}=1$, $g(\rho)^{2\phi(\delta)}=1$, $g(\rho)S(g,\chi)=1$.  Note that on some of the prime factors of $\rho$, the value of $g$ may have been chosen in an earlier round.

Then the corresponding completely multiplicative function $g$, with values in the complex unit circle, satisfies
\[
\lim_{x \to \infty} x^{-1} \sum_{n \le x} g(an+b) \overline{g}(An+B) = 1.
\]
Denoting $g(an+b) \overline{g}(An+B)$ by $z_n$ for convenience, we note that
\[
x^{-1} \sum_{n \le x} ( 1 - \Re z_n) \to 0, \quad x \to \infty.
\]
Here, if $z_n \ne 1$, $z_n^{2\phi(\delta)}=1$ guarantees that $1-\Re z_n \ge c_0 > 0$, the value of $c_0$ depending only upon $\delta$.

In particular, those $n$ for which $z_n \ne 1$ have asymptotic density zero.

Given an integer $n_0$ for which $a_1n_0+b_1 = w_1 r_1$ with $w_1$ comprised of primes dividing $\delta$, $r_1$ coprime to $\delta$, likewise $A_1 n_0 + B_1 = w_2 r_2$, then for a sufficiently large $k$ every integer of the class $n_0 \modd{\delta^k}$ will give rise to integers $a_1 n + b_1$, $A_1 n + B_1$ of the same form, with identical values of the $w_j$ and with their corresponding $r_j$ belonging to a respective fixed residue class $\modd{\delta}$.

Hence
\[
g\left( \frac{an+b}{An+B}\right) = g(\rho) \frac{g(w_1)}{g(w_2)} \frac{g(r_1)}{g(r_2)}
\]
has a constant value on a nonempty residue class.

Such a class has positive asymptotic density, hence contains a representative on which $g$ has the value 1.

Thus $g$ has the value 1 on the whole class.

In this way every $g\left( \frac{an+b}{An+B} \right)=1$, i.e. $g$ is a character on $G$.

This gives a characterisation of $G$ in terms of its dual group.

Consideration of the dual of the exact sequence $1 \to \Gamma \to \mathbb{Q}^* \to G \to 1$ shows the dual group of $\Gamma$ to be isomorphic to the quotient of denumerably many copies of the unit circle, one for each prime, by the dual of $G$.


\noindent \emph{Remarks}.  The choice of a principal character $\modd{\delta}$ enables the character $g$ identically 1 on $G$, so the conditions (i), (ii) and (iii) are consistent.

For a given character $\chi \modd{\delta}$ there can be at most one compatible set of values for $g$ on $\rho$ and the remaining torsion primes.  The ratio of two compatible values would yield a character on $G$ that is 1 on the primes not dividing $\delta$.  In view of Lemma \ref{elliott_kish_2015_groups_lem_01}, it would also be 1 on the primes not dividing $(a_1, A_1)$, hence on $\rho$.


\noindent \textbf{Determination of $G$; practical matters}.  The function $S(g,\chi)$ is given by an infinite series.  In this section it will be shown that it factorises and each of the finitely many factors is a polynomial in its associated variable $g(p)$.

Consider a typical sum
\[
\sum_{n \modd{\delta[d_1,d_2]}} \chi\left( \frac{a_1n + b_1}{d_1} \right) \overline{\chi} \left( \frac{A_1 n + B_1}{d_2} \right).
\]
Let $\delta = \prod_{j \le v} \ell_j^{\alpha_j}$, $\ell_j$ distinct primes, $d_1 = \prod_{j \le v} \ell_j^{\beta_j}$, $d_2 = \prod_{j \le v} \ell_j^{\gamma_j}$, where $\beta_j=0$, $\gamma_j = 0$ is possible, but $\min(\beta_j, \gamma_j)$ is bounded by the constraint $\ell_j^{\min(\beta_j, \gamma_j)} \mid \Delta_1$.

Consider the map
\[
n \modd{\delta[d_1, d_2]} \to \otimes \, n \modd{\ell_j^{\alpha_j + \max(\beta_j, \gamma_j)}}
\]
given by
\[
n \to \sum_{j \le v} u_j L_j,
\]
where $L_j \equiv 1 \modd{\ell_j^{\alpha_j + \max(\beta_j, \gamma_j)}}$, $L_j \equiv 0 \modd{\ell_r^{\alpha_r + \max(\beta_r, \gamma_r)}}$ if $1 \le r \le v$, $r \ne j$, (Chinese Remainder Theorem).

Typically $an + b \equiv 0 \modd{d_1}$ if and only if $an+b \equiv 0 \modd{\ell_j^{\beta_j}}$, i.e. $au_j + b \equiv 0 \modd{\ell_j^{\beta_j}}$, $1 \le j \le v$.  Similarly if $An + B \equiv 0 \modd{d_2}$.

Then
\[
\chi\left( \frac{a_1 n + b_1}{d_1} \right) = \prod_{j \le v} \chi_j \left( \frac{a_1 n + b_1}{d_1} \right),
\]
where $\chi_j$ is a Dirichlet character $\modd{\ell_j^{\alpha_j}}$.

In particular,
\[
\chi_j \left( \frac{a_1 n + b_1}{d_1} \right) = \chi_j \left( \frac{\sum_{r \le v} (a_1 u_r + b_1)L_r}{d_1} \right)
\]
since $\sum_{r \le v} L_r \equiv 1 \modd{\delta[d_1,d_2]}$.

Hence
\[
\chi_j \left( \frac{a_1 n + b_1}{d_1} \right) = \chi_j \left( \frac{a_1 u_j + b_1}{d_1} \right) = \chi_j \left( \frac{a_1 u _j +b_1}{\ell_j^{\beta_j}} \right) \overline{\chi}_j (d_{1j})
\]
where $d_{1j} = \ell_j^{-\beta_j} d_1$, $1 \le j \le v$; for $L_r/d_1$ is divisible by $\ell_j^{\alpha_j}$ if $r  \ne j$.

Define $\widehat{\chi}_j(\ell_j^{\beta_j})$ to be $\prod_{1 \le r \le v, r \ne j} \overline{\chi}_r(\ell_j^{\beta_j})$, $1 \le j \le v$.

The sum $S(g,\chi)$ becomes
\[
g(\rho) \prod_{j=1}^v \sum_{\substack{\beta_j \ge 0 \\ \gamma_j \ge 0}} \frac{g(\ell_j^{\beta_j}) \widehat{\chi}_j(\ell_j^{\beta_j}) \overline{g}(\ell_j^{\gamma_j})\overline{\widehat{\chi}}_j(\ell_j^{\gamma_j})}{ \ell_j^{\max(\beta_j, \gamma_j)}} \eta_j
\]
where
\[
\ell_j^{\alpha_j}\eta_j = \ell_j^{\alpha_j} \eta_j( \beta_j, \gamma_j) =  \sum_{u_j \modd{\ell_j^{\alpha_j + \max(\beta_j, \gamma_j)}}} \chi_j \left( \frac{a_1 u_j + b_1}{\ell_j^{\beta_j}} \right) \overline{\chi}_j \left( \frac{A_1 u_j + B_1}{\ell_j^{\gamma_j}} \right),
\]
it is understood that $a_1 u_j + b_1 \equiv 0 \modd{\ell_j^{\beta_j}}$, $A_1 u_j + B_1 \equiv 0 \modd{\ell_j^{\gamma_j}}$, and $\ell_j^{\min(\beta_j, \gamma_j)} \mid \Delta_1$, $1 \le j \le v$.

Since $g(\ell_j^\beta) = g(\ell_j)^\beta$, $\chi^\beta$ are periodic in $\beta$, period $2 \phi(\delta)$, a typical innersum becomes a polynomial in $g(\ell_j)\widehat{\chi}_j(\ell_j)$, of degree at most $2\phi(\delta)-1$, with coefficients that are linear forms in $2\phi(\delta)^{th}$ roots of unity (values of $\chi$) that in turn have coefficients that are essentially geometric progressions, indeed rational numbers.

Moreover, the constraint upon the size of $\min(\beta_j, \gamma_j)$ ensures that for all sufficiently large values of $\max(\beta_j, \gamma_j)$ either $\eta_j = 0$ or the character $\chi_j$ is principal.

The values of the $g(p)$, $p \nmid (a_1, A_1)$, together with that of $g(\rho)$, that fulfill conditions (ii) and (iii) may therefore be ascertained recursively.


\noindent \textbf{Further practical matters}.   Let $k$ be a positive integer.  The foregoing argument reduces the determination of the multiplicative group generated by the rationals $(an+b)/(An+B)$ with $n \ge k$ to a calculation in a polynomial ring over a cyclotomic extension of the rational field.

The following results may accelerate this process.

\begin{lemma} \label{elliott_kish_2015_groups_lem_extra}
Let the integers $a_j>0$, $b_j$ satisfy $\Delta_0 = a_1b_2 - a_2b_1 \ne 0$, and let the character $\chi \modd{p^m}$, with $p$ prime, have a constant value on the ratios $(a_1n + b_1)/(a_2n + b_2)$ for which $(a_j n + b_j, p)=1$, $j = 1, 2$.

Then $p$ divides $(a_1, a_2)\Delta_0$, or $p=3$ and the order of $\chi$ is at most 2, or $p=2$.
\end{lemma}

\noindent \emph{Proof}.  Suppose that $p$ does not divide $(a_1,a_2)\Delta_0$.  Then for each $w$ that satisfies $(w(a_1 - w a_2),p)=1$,
\[
n = \frac{wb_2 - b_1}{a_1 - wa_2} \qquad \text{gives} \qquad \frac{a_1n + b_1}{a_2 n + b_2} \equiv w \modd{p^m}.
\]
Note that $(a_2n + b_2)(a_1 - wa_2) \equiv \Delta_0$ and $(a_1n + b_1)(a_1 - wa_2) \equiv w \Delta_0 \modd{p^m}$, so that the ratio involving $n$ satisfies the requirement for $\chi$ to have a nonzero value.

Accordingly, $\chi(w)$ has a nonzero constant value on at least $p^m - 2p^{m-1}$ reduced residue classes $\modd{p^m}$.

If $\chi$ has order $t$, then it can have a fixed value on at most $t^{-1} \phi(p^m)$ reduced classes $\modd{p^m}$.

Hence $p^m - 2p^{m-1} \le t^{-1}  p^{m-1}(p-1)$; $t \le (p-1)/(p-2) < 2$ whenever $p >3$, i.e. $t=1$.  Otherwise $p=3$ and $t \le 2$, or $p=2$.

This completes the proof of the Lemma.

\noindent \emph{Example}.  The groups $\mathbb{Q}^*/\Gamma_k$ attached to the ratios $(3n+1)/(5n+2)$, $n \ge k$, introduced in the first author's volume on arithmetic functions and integer products \cite{Elliott1985}, were there shown, via homomorphisms into the positive reals, to be finite.

After Theorem \ref{elliott_kish_2015_groups_thm_01} and Lemma \ref{elliott_kish_2015_groups_lem_extra} each character $g$ on $\mathbb{Q}^*/\Gamma_k$ coincides, on the primes $p$ that do no divide 30, with a Dirichlet character $\chi \modd{3^3}$ that is of order at most 2.

The squares $\modd{3^3}$ are represented by $1, 4, 16, 25, 22, 10, 19, 13$ and $7$.  Since $(3 \cdot 4 + 1)/(5 \cdot 4 + 2) \equiv 13/2\cdot 11$, for $n \equiv 4 \modd{5\cdot 2^2 \cdot 3^3}$
\[
1 = g\left( \frac{3n+1}{5n+2} \right) = \chi(13) \overline{g}(2) \overline{\chi}(11)
\]
and, if $\chi$ is indeed quadratic, $g(2) = -1$.

However, for $n \equiv 9 \modd{5 \cdot 2^3 \cdot 3^5}$
\[
1 = g\left( \frac{3n+1}{5n+2} \right) = g(2^2) \chi(7) \overline{\chi}(47) = -1;
\]
a contradiction ensues.

The character $\chi$ is, in fact, principal and the argument shows $g(2)=1$.

Similarly, $n \equiv 2 \modd{2^3 \cdot 3^4 \cdot 5}$ gives $g(3)=1$ and $n \equiv 3 \modd{2^2 \cdot 3^3 \cdot 5^2}$ gives $g(5)=1$.

The groups $\mathbb{Q}^*/\Gamma_k$ are trivial.  Each positive rational $r$ has infinitely many representations
\[
r = \prod_j \left( \frac{3n_j+1}{5n_j+2} \right)^{\varepsilon_j}, \quad \varepsilon_j = \pm 1,
\]
with the $n_j$ as large as desired.  In the notation of the preface to the volume \cite{Elliott1985}, we may take $v=1$.

Rationalising the denominators establishes the

\noindent \textbf{Multiplicity Lemma.} \emph{Let the integers $a_j>0$, $b_j$, $j=1, 2$, satisfy $\Delta_0 = a_1b_2 - a_2 b_1 \ne 0$, and the integers $s$, $n$, $n'$, $(s, (a_2n+b_2)(a_2n' + b_2))=1$.}

\emph{Then}
\[
\frac{a_1n + b_1}{a_2n + b_2} \equiv \frac{a_1 n' + b_1}{a_2 n' + b_2} \modd{s}
\]
\emph{if and only if $\Delta_0(n-n') \equiv 0 \modd{s}$.}

\noindent \emph{Example}.  Consider the group $\mathbb{Q}_5^*/\Gamma$, where $\mathbb{Q}_5^*$ denotes the multiplicative positive rationals not divisible by 5 and $\Gamma$ its subgroup generated by all but finitely many fractions of the form $(5n+1)/(5n-1)$.  It was established in the first author's volume \cite{Elliott1985} that $\mathbb{Q}_5^*/\Gamma$ is finite.  Reduction $\modd{5}$ shows 2 not to belong to $\Gamma$, hence that $\mathbb{Q}_5^*/\Gamma$ is not trivial.

According to the main result of the present paper, on the primes that do not divide 10, each character $g$ on $\mathbb{Q}_5^*/\Gamma$ coincides with a product of Dirichlet characters, $\chi_1 \chi_2$, to moduli $2^4$ and $5^9$ respectively.

In particular,
\[
\chi_1\left( \frac{5^9 \cdot 2m + 1}{5^9 \cdot 2m - 1} \right) = g\left( \frac{5^9 \cdot 2m + 1}{5^9 \cdot 2m - 1} \right) \overline{\chi}_2(-1)
\]
has a constant value on all integers $m$.  After the Multiplicity Lemma the ratios $(5^9 \cdot 2m + 1)/(5^9 \cdot 2m - 1)$ determine a set of reduced residue classes $\modd{2^4}$, each class covered exactly 4 times.  If $\chi_1$ has order $t$, then $2^2 \le t^{-1}\phi(2^4)$, $t \le 2$.  The character $\chi_1$ is real.

Likewise, a consideration of $g((5^9 \cdot 2m + 1)/(5^9 \cdot 2m - 1))$ shows the order of $\chi_2$ not to exceed 4, i.e. to be 1, 2, or 4.

By the argument of the constraints section of the present paper, $g(2)$ assumes a value taken by $\chi_1 \chi_2$; $g^4$ is identically 1.

If $n \equiv 3 \modd{2^8 \cdot 5^8 \cdot 7}$, then
\[
1 = g\left( \frac{5n+1}{5n-1} \right) = g\left( \frac{2^4}{2 \cdot 7} \right) \chi_1 \chi_2 \left( \frac{2^{-4} (5n+1)}{2^{-1}(5n-1)} \right) = g(2)^3 \overline{g}(7);
\]
$g(2) = \overline{g}(7)$.

Choosing $n \equiv 1 \modd{2^6 \cdot 3 \cdot 5^8}$ we see that $g(3) \overline{g}(2) = 1$, i.e. $g(2) = g(3)$.

Choosing $n \equiv 1 \modd{2^7 \cdot 3^3 \cdot 5^8 \cdot 7}$ yields $g(2^3) g(7) \overline{g}(2) \overline{g}(3^3)=1$, i.e. $g(2) = g(7)$.

Hence $g(2) = g(3) = g(7)$, $g(3)^2 = g(2)^2=1$.  In particular, $1 = g(3)^2 = \chi_2(3)^2$.  Since $3$ is a primitive root $\modd{5^k}$ for $k \ge 1$, if $\chi_2$ is to have order 4, then $\chi_2(3) = \pm i$, $\chi_2(3)^2 = -1$.  The character $\chi_2$ is also real.  When defined, the indices of a prime $\modd{5^k}$ and $\modd{5}$, with respect to 3, have the same parity.  We may identify $\chi_2$ with the Legendre symbol $\modd{5}$.

Then $1 = g(7/3) = \chi_1(-3)$.  Taking powers of $-3$ we see that $\chi_1$ assumes the value 1 on each of the classes $13, 9, 5$ and $1 \modd{2^4}$.

Moreover, choosing $n \equiv 2 \modd{2^4 \cdot 3^2 \cdot 5^8 \cdot 11}$ we see that $1 = g(11/3^2) = \chi_1(11)$.  Multiplying by powers of $-3$, $\chi_1$ also assumes the value 1 on the classes $11, 15, 3$ and $7 \modd{2^4}$.  The character $\chi_1$ is principal.

To complete the characterisation, $g(2) = g(3) = \chi_2(3) = -1 = \chi_2(2)$.



It is interesting to derive this result by more directly following the calculation in the section on practical matters.

In the notation of the foregoing example we wish to calculate
\[
\eta( \beta, \gamma) = 2^{-4} \sum_{n \modd{2^{4+\max(\beta,\gamma)}}} \chi_1 \left( \frac{5n+1}{2^\beta} \right) \chi_1 \left( \frac{5n-1}{2^\gamma} \right)
\]
where $\chi_1$ is a character $\modd{2^4}$ of order 1 or 2.

Since $2 \| (5n+1)$, $2 \| (5n-1)$ cannot be simultaneously effected, $\eta(1,1)=0$.  Likewise $\eta(\beta, \gamma) = 0$ if $\max(\beta, \gamma) \ge 1$ with $\min(\beta, \gamma)=0$.  

The substitutions $(5n+1)/2^\beta = w \modd{2^4}$, $(5n-1)/2^\gamma = -t \modd{2^4}$, show that for $\max(\beta, \gamma) \ge 2$,
\[
\eta( \beta, \gamma) = 2^{-4} \sum_{w \modd{2^4}} \chi_1(w) \chi_1( -1 + 2^{\max(\beta, \gamma)-1}w).
\]
In particular, $\eta( \beta, \gamma)=0$ if $\max(\beta, \gamma) \ge 5$ and $\chi_1$ is nonprincipal.

Further, if $\chi_1$ is nonprincipal our basic condition (iii) in the characterisation of $\mathbb{Q}^*/\Gamma$ becomes
\[
\sum_{m=2}^4 2^{-m} \eta(m,1) (z^{m-1} + z^{-m+1}) + \eta(0,0) = \overline{\chi}_2(-1) S(g,\chi_1),
\]
with $z = g(2) \overline{\chi}_2(2)$.  Since, for a principal character $\modd{2^4}$, $\eta(\beta, \gamma)=1/2$ unless zero, the absolute value of the left hand side of this equation is at most $\sum_{m=2}^4 2^{-m} + 2^{-1} = 2^{-4}(15) < 1$, which is untenable.

Therefore the character $\chi_1$ is principal.

The basic condition (iii) simplifies to $(2-z)\overline{\chi}_2(-1) = 1$.  Hence $\chi_2(-1) = 1$, so that $\chi_2$ is not quartic, but may be quadratic.  Moreover, $z=1$, $g(2) = \chi_2(2)$.


The group $\mathbb{Q}^*/\Gamma$ generated by the ratios $(5n+1)/(5n-1)$ has the single free generator, $5$, and a torsion group of order 2 determined by its dual through the quadratic Dirichlet character $\modd{5}$.

There are infinitely many representations
\[
57^2 = \prod_j \left( \frac{5n_j+1}{5n_j-1} \right)^{\varepsilon_j}, \quad \varepsilon_j = \pm 1,
\]
but no such representation is available to 57 itself.


\noindent \emph{Remarks}.  A short elementary proof that a complex-valued multiplicative function constant on all sufficiently large members of a progression $an+b$, $a>0$, coincides with a Dirichlet character $\modd{a}$, on the integers prime to $a$, may be found as Lemma 19.3 in \cite{Elliott1985}, pp. 334--335.

Let $\mathbb{Q}_2 = \mathbb{Q}^* \oplus \mathbb{Q}^*$ be the direct sum of two copies of the multiplicative positive rationals, $\Gamma_2$ its subgroup generated by the pairs $(an+b) \oplus (An+B)$, $n > k$.  Simultaneous representations of the form
\[
r_1 = \prod_{j=1}^m (an_j + b)^{\varepsilon_j}, \qquad r_2 =\prod_{j=1}^m (An_j + B)^{\varepsilon_j}, \quad \varepsilon_j = \pm 1,
\]
may be studied through the offices of the quotient group $\mathbb{Q}_2/\Gamma_2$.  A typical character on that group amounts to a pair of completely multiplicative functions $g_1$, $g_2$, with values in the complex unit circle, and that satisfy
\[
g_1(an+b) g_2(An+B) = 1, \quad n > k.
\]

We may reduce this two-dimensional problem to a one-dimensional problem by means of the following argument, given in an equivalent form in \cite{Elliott1985}, Chapter 19.  For ease of notation we shall assume, as we clearly may, that $b>0$, $B>0$.

Replacing $n$ by $bBn$,
\[
g_1(aBn + 1) g_2(Abn+1) = c_1 \ne 0.
\]

Replacing $n$ by $(aB+1)n+1$,
\[
g_1(aBn+1)g_2(Ab[(aB+1)n+1]+1) = c_2 \ne 0.
\]

Eliminating between these relations,
\[
g_2(Ab(aB+1)n+1)\overline{g}_2(Abn+1) = c_3 \ne 0,
\]
to which we can apply Theorem \ref{elliott_kish_2015_groups_thm_01}.  The functions $g_2$ and $g_1$ are essentially Dirichlet characters and we may follow the treatment for a single function.


There is an alternative procedure in Steps \ref{elliott_kish_2015_groups_step_02} and \ref{elliott_kish_2015_groups_step_03} that may more readily generalise to higher dimensional problems.  Once Step \ref{elliott_kish_2015_groups_step_02} guarantees the existence of a real $\alpha$ for which the series $\sum p^{-1} ( 1 - \Re g(p)^k p^{-i\alpha})$ converges, Lemma \ref{elliott_kish_2015_groups_lem_C} is applied to the completely multiplicative function $h$ defined by $p \to h(p)=g(p)^k p^{-i\alpha}$.  We arrive at an asymptotic estimate
\[
\lim_{x \to \infty} x^{-1} \sum_{n \le x} h(an+b) \overline{h}(An+B) = \prod_p w_p
\]
where the $w_p$ may now depend upon the parameter $\alpha$.  Following the argument of Step \ref{elliott_kish_2015_groups_step_03}, $g(p)^k = p^{i\alpha}$ on all but finitely many primes.

Bearing in mind the two cases considered in Step \ref{elliott_kish_2015_groups_step_01}, typically, on a suitable residue class $(g(an+b)\overline{g}(An+B))^k$ will coincide with $((an+b)/(An+B))^{i\alpha}$.

A simple analogue of Theorem \ref{elliott_kish_2015_groups_thm_01} then suffices to guarantee that $\alpha=0$:

\begin{lemma} \label{elliott_kish_2015_groups_lem_extra6}
If integers $u_j >0$, $v_j$, $j=1, 2$, $\Delta = u_1v_2 - u_2v_1 \ne 0$, satisfy
\[
\left( \frac{u_1n + v_1}{u_2n + v_2} \right)^{i\alpha} = c
\]
for all $n$ sufficiently large, then $\alpha=0$ and $c=1$.
\end{lemma}

\noindent \emph{Proof of Lemma \ref{elliott_kish_2015_groups_lem_extra6}}.  Since $(u_j n + v_j)^{i\alpha} = (u_j n)^{i\alpha} ( 1 + i\alpha v_j (u_j n)^{-1} + O(n^{-2}))$, $j = 1, 2$,
\[
c = \left( \frac{u_1}{u_2} \right)^{i\alpha} \left( 1 + \frac{i\alpha \Delta}{u_1u_2 n} + O\left( \frac{1}{n^2} \right) \right), \quad n \to \infty.
\]
This forces $(u_1/u_2)^{i\alpha}=c$, $\alpha=0$, $c=1$ in turn.

We may now continue as before until reaching the section on constraints.  To once again avoid applying Step \ref{elliott_kish_2015_groups_step_01}, choice of an appropriate residue class for $n$ will, for example, for a given prime $p$ arrange infinitely many representations
\[
1 = g\left( \frac{an+b}{An+B} \right) = g(p) \chi\left( \frac{an+b}{p(An+b)} \right)
\]
with $p \| (an+b)$, $(p^{-1}(an+b)(An+B),\delta)=1$, where $\chi$ is a Dirichlet character $\modd{\delta}$.

In particular, $g(p)$ will be a value of the character $\chi$, and $g(p)^k=1$ with $k$ the order of $\chi$.

Altogether, this variant argument obviates appeal to the Fourier analysis, involving estimates for Kloosterman sums, that is given in \cite{Elliott1985}, Chapters 2 and 4.  It does not deliver, for the moment at least, the upper bound on the number of terms sufficient to represent group theoretically a typical positive rational, $r$, or the underlying recursive argument by which such a bound is obtained there.  Moreover, the attendant inequalities on the additive functions in \cite{Elliott1985}, Chapter 10, are strongly localised.


\noindent \emph{Example}.  Let the integers $a>0$, $A>0$, $b$, $B$ satisfy $\Delta = aB - Ab \ne 0$, $(a,b)=1=(A,B)$.  In order to avoid complications involving the prime 2, we assume that if $aA$ is odd, then $b$ and $B$ have the same parity.  Then there is an integer $n_0$ for which $(an_0 + b)(An_0 + B)$ is coprime to $\delta = 6aA\Delta$.

Let $\mathbb{Q}_\delta^*$ be the multiplicative group of positive rationals coprime to $\delta$, and $\Gamma_\delta$ its subgroup generated by all but finitely many fractions $(an+b)/(An+B)$, $n>0$, where $n$ belongs to the residue class $n_0 \modd{\delta}$.  This ensures that every member of $\Gamma_\delta$ belongs to $\mathbb{Q}_\delta^*$.

Set $G_\delta = \mathbb{Q}_\delta^*/\Gamma_\delta$.

A typical character $g$ on $G_\delta$ satisfies $g\left( \frac{an+b}{An+B} \right)=1$ on all but finitely many positive representatives of the class $n_0 \modd{\delta}$.  By Theorem \ref{elliott_kish_2015_groups_thm_01} there is a Dirichlet character $\chi\modd{\delta^4}$ for which $g(p) = \chi(p)$ whenever $(p,\delta^4)=1$.  Since this comprises precisely the primes in $\mathbb{Q}_\delta^*$, $g$ is largely determined.

Moreover, any character $\modd{\delta^4}$ that is identically 1 on the fractions $(an+b)/(An+B)$, $n \equiv n_0 \modd{\delta}$, induces a character on $G_\delta$.

The group $G_\delta$ is thus isomorphic to the group of such Dirichlet characters.

\noindent \emph{Closing Remarks (by the first author)}.  In anticipation of the validity of an appropriate version of the conjecture that correlations of multiplicative functions with values in the complex unit circle could only satisfy
\[
\limsup_{x \to \infty} x^{-1} \left| \sum_{n \le x} g_1 (an+b) g_2(An+B) \right| > 0, \quad \Delta \ne 0,
\]
if there are reals $\tau_j$ and Dirichlet characters $\chi_j$ so that the series
\[
\sum p^{-1} (1- \Re g_j(p) \chi_j(p) p^{i\tau_j}), \quad j = 1, 2,
\]
taken over the prime numbers, converge (c.f. \cite{MR1042765} Conjecture II, p. 68, \emph{convolutions} in the title of that paper should have read \emph{correlations}; modified slightly in \cite{elliott1994correlationmemoir}  Conjecture II, p. 4), towards which the result of Terence Tao is a decisive step, Jonathan Kish and I were already in possession of a detailed version of Step \ref{elliott_kish_2015_groups_step_04} of the present account by mid 2010 \underline{\hspace{.2in}} putting the cart before the horse is sometimes helpful.

The conjecture has a root in the Probabilistic Theory of Numbers.  In the early nineteen seventies, when studying asymptotic behaviour of additive arithmetic functions with unbounded renormalisations:
\[
\mathop{ \phantom{1} [x]^{-1} \sum 1 \phantom{[x]^{-1}}}_{\substack{n \le x \\ f(n) - \alpha(x) \le z \beta{x}}}, \quad x \, \text{(real)}\, \to \infty,
\]
I felt that the Erd\H os-Kac realisation of an additive function as a sum of independent random variables, in general not valid, might be restored provided a suitable (moving) obstruction were removed from $f$.

Where to find such an obstruction?  In the event, c.f. \cite{elliott1975lawlargenos}, \cite{elliott1976generalasymptotic}, perturbation of the underlying probability models, somewhat in the style of the perturbation of planetary orbits within the circle of ideas of the Hamilton-Jacobi equation, led to a satisfactory outcome, the renormalising functions $\alpha(x)$, $\beta(x)$ classified according to their behaviour under the group of transformations $x \to x^y$, $y>0$, fixed.

The same point of view could be applied to frequencies involving the sums $f_1(an+b) + f_2(An+B)$ of two (or more) possibly distinct additive functions, $aB - Ab \ne 0$, c.f. \cite{elliott1979sumsdifferencesadditivemeansquare}.

The respective characteristic functions of the corresponding frequencies have the form
\[
[x]^{-1} \sum_{n \le x} g(n), \qquad [x]^{-1} \sum_{n \le x} g_1(an+b) g_2(An+B),
\]
where the functions $g(n) = \exp(it f(n))$, $g_j(n) = \exp(it f_j(n))$, $t$, $t_j$ real,  are multiplicative.

Besides the many experiences in Probabilistic Number Theory, there is the further experience of the Hardy-Ramanujan-Littlewood Circle Method.  One may view that method as resting on the approximation of continuous characters $\alpha \to \exp(2\pi i k \alpha)$, $k = 0 \pm 1, \pm 2, \dots$, on the group $\mathbb{R}/\mathbb{Z}$ dual to the additive group of integers, by the discrete characters $a/q \to \exp(2\pi i ka/q)$, (usually effected with $(a,q)=1$) on the additive (reduced) group $\mathbb{Z} \modd{q\mathbb{Z}}$.

For this and other reasons the Dirichlet characters suggested themselves as obstructions in the dual group of the positive rationals under multiplication, the direct product of denumerably many copies of $\mathbb{R}/\mathbb{Z}$, and in which, in a certain sense, they are dense, c.f. \cite{Elliott1997} Chapter 12, Exercise 7.  The role of reduced rationals $a/q$ is then played by primitive characters.

Whilst the Stone-Weierstrass theorem might offer other approximating functions, the connections between special functions and group representations also suggests the application of associated group characters.

Conjecture III of \cite{elliott1994correlationmemoir}, p. 65, that there is a positive absolute constant $c$ such that
\[
x^{-1} \sum_{n \le x} g(n) h(n+1) \ll \left( T^{-1} + \exp\left( -\min_\chi \min_{|\tau| \le T} \sum_{p \le x} p^{-1} (1 - \Re g(p) \chi(p) p^{i\tau} )\right)  \right)^c
\]
uniformly for $x \ge 2$, and multiplicative functions $g$, $h$ with values in the complex unit disc, modified by the requirement that the characters $\chi$ be to moduli not exceeding $T$, as noted in the author's Cambridge Tract \cite{Elliott1997}, Chapter 34, p. 315, may well be applied to $g(an+b)h(An+B)$, $aB - Ab \ne 0$.

It is clear that some restriction must be made upon the size of the defining moduli of the characters $\chi$, since the application of Kronecker's theorem that shows the $\chi$ to be dense in the dual of $\mathbb{Q}^*$ for a given value of $x$ shows the minimum over $\chi$ in Conjecture III to be arbitrarily close to zero.

Although the ultimate aim was for a fixed obstruction, experience in Probabilistic Number Theory showed that the variable $\tau$ might be allowed to float far past $x$ in size yet be retrieved, c.f. \cite{elliott1988localizederdoswintner}, \cite{elliott2010valuedistribution}.  Moreover, c.f. \cite{elliottMFoAP1}, \cite{elliottMFoAP6}, I was aware that for considerable ranges of $\tau$ and the modulus $D$, at most one generalised character $n \to \chi(n)n^{i\tau}$, $\chi\modd{D}$, could be near to a given multiplicative function.  To this extent the obstruction would be isolated (just as it is in the case of a single renormalised additive function).

Ultimately, one would expect the study of quotient groups of $\mathbb{Q}^*$ to embrace integration over appropriate subgroups of the dual group of $\mathbb{Q}^*$.

An extensive survey of problems and results attached to many dimensional product representations of rationals by the values of polynomials on the integers or on the primes, together with a discussion of their attendant groups and $\mathbb{Q}^*$-character sums, may be found in \cite{elliott2002productrepresentations}.

\bibliographystyle{amsplain}
\bibliography{MathBib}

\end{document}